\input amstex 
\documentstyle{amsppt}
\input bull-ppt
\keyedby{bull425e/PAZ}

%%%%%%% Author macros %%%%%%%%%%%%%%%%%%%%%%
	\def\integer{{\Bbb Z}}
	\def\real{{\Bbb R}}
	\def\iitem{\itemitem}
%%%%%%% End of Author macros %%%%%%%%%%%%%%%%%%%%%%

\topmatter
\cvol{29}
\cvolyear{1993}
\cmonth{October}
\cyear{1993}
\cvolno{2}
\cpgs{213-217}
%\ratitle
\title Symmetry of Tilings of the Plane\endtitle
\author Charles Radin \endauthor
%\shortauthor{}
%\shorttitle{}
\address Mathematics Department,
University of Texas,
Austin, Texas 78712\endaddress
\ml radin\@math.utexas.edu \endml
%\cu \endcu
\date February 12, 1993\enddate
\subjclass Primary 52C20, 58F11, 47A35\endsubjclass
\thanks Research supported in part by Texas ARP Grant 
003658-113\endthanks
\abstract We discuss two new results on tilings of the 
plane.
In the first, we give sufficient conditions for the tilings 
associated with an inflation rule to be uniquely ergodic
under translations, the conditions holding for the pinwheel
inflation rule. In the second result we prove there are
matching rules for the pinwheel inflation rule, making the
system the first known to have complete rotational 
symmetry.\endabstract
\endtopmatter

\document

We consider tilings of the Euclidean plane, $E^2$, by 
(orientation-preserving) congruent copies
of a fixed finite set of prototiles. Prototiles are 
topological disks
in the plane satisfying some mild restrictions on their 
shapes, as
detailed below. Congruent copies of prototiles are called 
tiles, and a
tiling is simply an unordered collection of tiles whose 
union is the
plane and in which each pair of tiles has disjoint 
interiors.

We are concerned here with two constructions associated 
with a fixed
finite set $S=\{P_j\}$ of prototiles, the most important 
of which is
the set $X(S)$ of all tilings by tiles from $S$. In 
particular, we are
interested in understanding the purest cases, in which all 
the
tilings in $X(S)$ are ``essentially the same''; we will 
define this
precisely further on. Two examples are exhibited in 
Figures 1 and 2
on the next page, both with
two prototiles; in Figure 1, $S = S_K$ produces only a
checkerboard-like tiling (and all congruences), and in 
Figure 2,
$S=S_P$ produces the well-known tilings of Penrose [3, 4, 
6].

Tilings like those of Penrose are not usually invariant 
under any congruence
of the plane (other than the identity), so to analyze 
their symmetries
we introduce some elementary ergodic theory and another 
basic
construction which can sometimes be associated with a 
prototile set
$S$, the set $X^F(S)$ of tilings defined by an ``inflation 
function''
$F$.  An ``inflation rule'' for $S$, if it exists, 
consists of a
dilation $D_F$ of $E^2$ by some factor $\lambda_F < 1$ and 
a finite set
$\{C_{jk}\}$ of congruences of $E^2$, such that for each 
$P_j \in S$ we have
$$P_j=\bigcup_kC_{jk}D_FP_{n_k} \tag 1$$
\noindent
where the elements of each union have pairwise disjoint 
interiors.
The inflation function $F$ associated with such a rule is 
defined on
tiles (and then sets of tiles), with sets of tiles as 
values, as
follows. If the tile $P$ is ``of tile-type $j$'', that is, 
$P = CP_j$
where $P_j\in S$ and $C$ is a congruence, then 
$$F: P\longrightarrow F(P) 
\equiv\{D_F^{-1}CC_{jk}D_FP_{n_k}\}.
\tag2$$
\topspace{6pc}\caption{Figure 1\qquad\qquad\qquad
\qquad\qquad\qquad Figure 2}

\noindent
(Intuitively, $P$ can be replaced by a set of ``small-size 
tiles'' 
by (1), which are
then expanded in (2) to original size by the inverse
of the dilation. This process can obviously be applied to 
any collection of
tiles---for example, a tiling---and can thus be iterated.)

The tilings $X^F(S)$ associated with the inflation 
function $F$ are
then defined as those tilings $T$ such that each finite 
subcollection
of $T$ is congruent to a subcollection of a set of tiles 
of the form
$F^r(P)$ for some prototile $P$ and integer $r\geq 1$. 

The construction $X^F(S)$ is mainly of interest when $F$ 
is a
homeomorphism on it, since then $F$ defines a natural 
representation
for the dilation $D_F$ as a map of $X^F(S)$ onto itself. 
(It is ``natural'' in that it extends the representation 
of the
congruences to a larger subgroup of the conformal group.)
One can then
consider ``symmetry'' with respect to this hierarchical 
action, as we
shall do.  It is noteworthy that dilational (and rotational)
symmetry is manifested not by the invariance of tilings 
themselves
but rather of measures on tilings; that is, the action of 
the dilation
(and rotations) is lifted to the set of (translation 
invariant Borel
probability) measures on the tilings, and invariance is 
sought in this
set. This is one reason for using the machinery of ergodic 
theory in
the analysis of tilings.

We note the following examples of inflation rules.

\iitem{(i)\phantom{ii}} The square inflation rule $S$, 
given in Figure 3---$S_S$ has
one element;

\iitem{(ii)\phantom{i}} The Robinson inflation rule $R$, 
given in Figure 4---$S_R$ has
two elements;

\iitem{(iii)} Conway's ``pinwheel'' inflation rule $C$, 
given in Figure 
5---$S_C$ has two elements.

(The inflation rule is one-to-one---and in fact a
homeomorphism---on the associated space of tilings for the 
Robinson
and pinwheel rules but not the square rule, for which it is
four-to-one.
Also, we really should
distinguish between left-
and right-handed triangles in Figure 4.) 
Each $X(S)$ (and $X^F(S)$) carries a natural metric
structure (indicated below) and defines a dynamical system 
with $\real^2$
action, where $\real^2$ acts by simultaneous translation 
of the tiles in a
tiling.  Using this language, we note the unobvious fact 
[4] that the
dynamical systems associated with $X^R(S_R)$ and $X(S_P)$ 
are
topologically conjugate: there is a homeomorphism between 
$X^R(S_R)$
and $X(S_P)$ which\ intertwines the
translations. It is immediate that
those associated with $X^S(S_S)$ and $X(S_K)$ are 
topologically
conjugate.

\fighere{7pc}\caption{Figure 3}

\topspace{14pc}\caption{Figure 4}

We now outline our results. Our first result gives 
conditions
sufficient for the dynamical system associated with some 
$X^F(S)$ to
be uniquely ergodic, the conditions being satisfied, for 
example, by
$X^C(S_C)$. Unique ergodicity of a dynamical system, 
namely, the
property that there is one and only one Borel probability 
measure on
the space invariant under the group action---the space 
here being
$X^F(S)$ and the group being $\real^2$---is useful because 
it means
the dynamical system has a natural measure associated with 
it. (From
Birkhoff's pointwise ergodic theorem [1], this is 
equivalent to having
all the tilings being statistically identical [6].) The 
Penrose (and
therefore Robinson) systems are not quite uniquely 
ergodic. However,
each ergodic measure is invariant under rotation by 
$2\pi/10$ about
any point of $E^2$, and the different measures are merely 
rotations of one
another, by an angle in $(0,2\pi/10)$. So modulo this 
rotation there
is a unique invariant measure for these systems.  
Furthermore, the
ergodic measures for the Robinson system are each 
invariant under the
dilation associated with Figure~4.

Our second result is the construction of a finite 
prototile set
$\widetilde S$ such that the dynamical system associated 
with $X(\widetilde
S)$ is uniquely ergodic, and metrically conjugate with 
that associated
with $X^C(S_C)$.  (Metric conjugacy---the existence of a 
measurable
bijection modulo measure zero, which intertwines the 
dynamics---is
somewhat weaker than topological conjugacy.) In other words,
$X^C(S_C)$ plays a role for $X
(\widetilde S)$ similar 
to the one
$X^R(S_R)$ plays for $X(S_P)$; in each case the 
former\vadjust{\fighere{6pc}\caption{Figure 5}}\ defines 
for the
latter a hierarchical symmetry which, for example, is 
manifested in a
symmetry in the spectrum (discussed below) of the latter 
in each
ergodic component. The pinwheel and Penrose systems differ
significantly in their rotational symmetry: the Penrose 
system has the
tenfold rotational symmetry of each of its ergodic 
components, while
the pinwheel system, being uniquely ergodic, has full 
rotational
symmetry.  (The actions on a space of tilings of 
congruences of
$E^2$, and of a dilation $D_F$---assuming $F$ is a 
homeomorphism on
$X^F(S)$, of course---can be lifted to actions on the set 
of invariant
measures of the dynamical system; so if a dynamical system 
is uniquely
ergodic, its measure must be invariant under such 
actions.) The
symmetries of the spectrum of tiling systems---in 
particular, the
Penrose tilings---have been of major importance in their 
connection
with theories of the structure matter, for example, 
quasicrystals [6].

We now state our assumptions and results more fully.
We assume that the prototiles in $S$
satisfy the following conditions, besides being homeomorphs
of the closed unit disk:

\roster
\item"$\bullet$" \<$X(S)$ is nonempty; 
\item "$\bullet$" in the tilings being considered, the 
boundary of 
each tile can be covered by tiles in only finitely many 
ways, up to
congruence;
\item "$\bullet$" each prototile $P$ has small 
surface-to-volume ratio in 
the sense that:
$${{\operatorname{area}\{x\in tP\ :\ \Vert x-y\Vert \le 1\ 
\text{ for some }y\in
\partial(tP)\}}\over {\operatorname{area} 
\{tP\}}}\longrightarrow 0 
\tag 3$$
as the expansion factor $t \to \infty$.
\endroster

A metrizable topology is put on the space $X(S)$ or 
$X^F(S)$ of
tilings with the following neighborhood basis
$\{N_T(\epsilon)\,:\,\epsilon > 0\}$ of each tiling $T$: 
$N_T(\epsilon)$
consists of all tilings $T'$ such that within the circle 
in the plane
centered at the origin and of radius $1/\epsilon$, 
wherever $T$ has a
tile, so does $T'$, within distance $\epsilon$ in the 
Hausdorff metric
on compact sets. It is known that $X(S)$ and $X^F(S)$ are 
compact in
this topology and that the actions of translations 
$\real^2$ on
$X(S)$ and $X^F(S)$ are continuous [9]. Given the above 
then, our
results are the following.

\proclaim{Theorem 1 \rm\cite{8}}
Assume given some $X^F(S)$ as above, and assume
there is some $r\ge 1$ such that\/\RM:
\vskip1pt
\roster
\item"(a)"
for each prototile $P$, $F^r(P)$ contains tiles of every 
tile-type\/\RM;
\vskip1pt
\item"(b)" for some prototile $P,$ $F^r(P)$ contains two 
tiles $P',\ P''$ of the 
same tile-type, whose relative rotation is irrational with 
respect to $\pi$.
\endroster
\vskip1pt
\noindent Then $(X^F(S),\real^2)$ is uniquely ergodic.
\endproclaim

\proclaim{Theorem 2 \rm\cite{7}}
There is a finite set $\widetilde S$ of prototiles
such that $(X(\widetilde S), \real^2)$ is uniquely ergodic 
and metrically
conjugate to $(X^C(S_C), \real^2)$.
\endproclaim

We note that $X^C(S_C)$ satisfies the hypotheses of 
Theorem 1 (with
$r=2)$ and close with a few comments on these results of a 
more
technical nature.

Much of the symbolic (substitution) dynamics of hierarchical
structures makes use of a square matrix $A$, of which 
$A_{jk}$ is the
number of times tile-type $j$ is associated, by the 
inflation function
$F$, with tile-type $k$. To prove Theorem 1, this matrix 
is generalized to
a family $A[m]$ of matrices for which $A[m]_{jk}$ is the sum
$\sum_ne^{ima_n(j,k)}$, over all tiles of type $j$ 
contained in
$F(P_k)$, where $a_n(j,k)$ is the angle of rotation of the 
tile
compared to the defining prototile of its type.

This generalization allows us to keep track of rotational 
information,
and we prove unique ergodicity by using Weyl's criterion 
on uniform
distribution as it is used in elementary treatments of 
rotations of
the circle [1]. The natural way in which the important 
matrix $A$
generalizes is evidence to us that much of symbolic 
substitution
dynamics can be generalized to tilings $X^F(S)$.

The import of Theorem 2 is also, in part, its relation to 
symbolic
dynamics. The systems $X(S)$ are natural generalizations 
of symbolic
systems of finite type with $\integer^2$ action; for 
example, they play
roughly the same role for statistical mechanics in 
Euclidean spaces as
systems of finite type do for lattice gas statistical 
mechanics [6]. Of
particular importance, for example, in theories of 
material structure
are the spectra of such dynamical systems. To see why the 
symmetries
of the dynamical system are relevant for this, note that 
if a rotation or
dilation $W$ preserves the ergodic probability measure $m$ 
of the
dynamical system, then, as shown in [8], the spectral 
projection
$E_{\Delta}$ of the translations, associated to a set 
$\Delta\subset \real^2$, is unitarily
equivalent, in the Hilbert space defined by $m$, to 
$E_{W(\Delta )}$.
(See [2] for explicit connections between dynamics spectra 
and X-ray
spectra of scatterers.) 

Another aspect of Theorem 2 is its relation to the work of 
S.\ Mozes,
who proved in [5] that a rather general class of symbolic 
substitution
dynamical systems with $\integer^2$ action are of finite 
type; our example
is, we hope, the first step in a parallel theorem for 
tilings of the
plane.

\enddocument